\documentclass[12pt,leqno,fleqn]{amsart}  
\usepackage{tikz}
\usepackage{amsmath,amstext,amsthm,amssymb,amsxtra}
\usepackage[top=1.5in, bottom=1.5in, left=1.25in, right=1.25in]	{geometry}
\usepackage{txfonts} 
\usepackage[T1]{fontenc}
\usepackage{lmodern}

\usepackage{mathtools}
\mathtoolsset{showonlyrefs,showmanualtags}

\usepackage{hyperref} 
\hypersetup{
    colorlinks=true,       
    linkcolor=blue,          
    citecolor=magenta,        
    filecolor=magenta,      
    urlcolor=cyan           
}

\usepackage[msc-links]{amsrefs}

\theoremstyle{plain} 
\newtheorem{lemma}{Lemma}[section]
\newtheorem{proposition}{Proposition}[section] 
\newtheorem{theorem}{Theorem} [section]
\newtheorem{corollary}{Corollary} [section]
\newtheorem*{question}{Question}
\newtheorem*{remark}{Remark}
\newtheorem{priorTheorem}{Theorem}

\def\ZM{\ensuremath{\mathcal M}}

\def\ZA{\ensuremath{\mathcal A}}

\def\ZI{\ensuremath{\mathbf 1}}
\def\ZN{\ensuremath{\mathbb N}}

\def\ZP{\ensuremath{\mathcal P}}
\def\ZH{\ensuremath{\mathcal X}}
\def\kappa{\ensuremath{\mathcal K}}
\def\ZR{\ensuremath{\mathbb R}}

\def\ZF{\ensuremath{\mathcal F}}
\def\ZP{\ensuremath{\mathcal P}}

\def\ZH{\ensuremath{\mathcal H}}

\theoremstyle{plain}
\newtheorem{definition}[equation]{Definition} 

\numberwithin{equation}{section}

\newcommand {\e }[1]{\eqref{#1}}
\newcommand {\lem }[1]{Lemma \ref{#1}}
\newcommand {\cor }[1]{Corollary \ref{#1}}
\newcommand {\pro }[1]{Proposition \ref{#1}}
\newcommand {\trm }[1]{Theorem \ref{#1}}

\title[] {On systems of non-overlapping Haar polynomials }
\author{Grigori A. Karagulyan}
\address{Faculty of Mathematics and Mechanics, Yerevan State
University, Alex Manoogian, 1, 0025, Yerevan, Armenia} 
\email{g.karagulyan@ysu.am}

\thanks{Research was supported by the Science Committee of Armenia, grant 18T-1A081 }

	
\begin{document}
\begin{abstract}
We prove that $\log n$ is an almost everywhere convergence Weyl multiplier for the orthonormal systems of non-overlapping Haar polynomials. Moreover, it is done for the general systems of martingale difference polynomials.
\end{abstract}

	\maketitle  
\section{Introduction}

The following two theorems are well-known in Fourier Analysis.
\begin{priorTheorem}[Menshov-Rademacher, \cite{Men}, \cite{Rad}, see also \cite{KaSa}]\label{MR}
If $\{\phi_k:\, \,k=1,2,\ldots,n\}\subset L^2(0,1)$ is an orthogonal system, then
\begin{equation}\label{1}
\left\|\max_{1\le m \le n}\left|\sum_{k=1}^m\phi_k\right|\,\right\|_2\le c\cdot\log n \left\|\sum_{k=1}^n\phi_k\right\|_2,
\end{equation}
where $c>0$ is an absolute constant.
\end{priorTheorem}
\begin{priorTheorem}[Menshov, \cite{Men}]\label{M}
 For any $n\in\ZN$ there exists an orthogonal system $\phi_k$, $k=1,2,\ldots,n$, such that
	\begin{equation}
	\left\|\max_{1\le m \le n}\left|\sum_{k=1}^m\phi_k\right|\,\right\|_2\ge c\cdot\log n \left\|\sum_{k=1}^n\phi_k\right\|_2,
	\end{equation}
	for an absolute constant $c>0$.
\end{priorTheorem}
Let $\label{Phi}
\Phi=\{\phi_k(x),\,k=1,2,\ldots\}\subset L^2(0,1) $ be an infinite orthogonal system of functions. Denote by $\ZP_n(\Phi)$ the family of all \textit{monotonic} sequences of $\Phi$-polynomials
\begin{equation*}
p_k(x)=\sum_{j\in G_k}c_j\phi_j(x),\quad k=1,2,\ldots, n,
\end{equation*}
where $G_1\subset G_2\subset \ldots \subset G_n\subset \ZN$ and $\sum_{j\in G_n}c_j^2\neq 0$. Define 
\begin{equation}
\kappa_n(\Phi)=\sup_{\{p_k\}\in \ZP_n(\Phi)}\frac{\left\|\max_{1\le m\le n}\left|p_m\right|\right\|_2}{\left\|p_n\right\|_2}.
\end{equation} 
From \trm{MR} it follows that $\kappa_n(\Phi)\le c\cdot \log n$  for every orthogonal system $\Phi$, where $c$ is an absolute constant.
On the other hand, applying \trm{M}, one can also construct an infinite orthogonal system with the lower bound $\kappa_n(\Phi)\ge c\cdot \log n$, $n=1,2,\ldots $. Thus we conclude, in general, the logarithmic bound of $\kappa_n(\Phi)$  is optimal. We will see below that from results of Nikishin-Ulyanov \cite{NiUl} and Olevskii \cite{Ole}  it follows that $\kappa_n(\Phi)\gtrsim \sqrt{\log n}$ for any complete orthonormal system $\Phi$. 

In this paper we found the sharp rate of the growth of $\kappa_n\sim \sqrt{\log n}$ for the generalized Haar systems. The classical Haar system case of the result is also new and interesting. The upper bound $\kappa_n\lesssim \sqrt{\log n}$ is proved for the general systems of martingale type. 

To state the main results recall few standard notations. The relation $a\lesssim b$ ($a\gtrsim b$) will stand for the inequality $a\le c\cdot b$ ($a\ge c\cdot b$), where $c>0$ is an absolute constant. Given two sequences of positive numbers $a_n,b_n>0$, we write $a_n\sim b_n$ if we have $c_1\cdot a_n\le b_n\le c_2\cdot a_n$, $n=1,2,\ldots$ for some constants $c_1, c_2>0$.  Throughout the paper, the base of $\log$ is equal $2$. 
\begin{theorem}\label{T1}
	If $\Phi$ is a martingale difference, then  $\kappa_n(\Phi)\lesssim \sqrt{\log n}$.
\end{theorem}
\begin{theorem}\label{T2}
	For any generalized Haar system $\ZH$ we have the relation 
	\begin{equation}\label{a2}
	\kappa_n(\ZH)\sim \sqrt{\log n}.
	\end{equation}
\end{theorem}
 In the class of all martingale differences the upper bound in \trm{T1} is optimal that readily follows from \trm{T2}.  One can easily see that for the Rademacher system we have $\kappa_n\sim 1$. So relation \e{a2} can not be extended for general martingale differences. Such estimates of $\kappa_n(\Phi)$ characterize Weyl multipliers of a given orthonormal system $\Phi$.  Recall some well-known definitions in the theory of orthogonal series (see \cite{KaSa}).
\begin{definition}
	Let $\Phi=\{\phi_n:\, n=1,2,\ldots\}$ be an orthonormal system ($\|\phi_n\|_2=1$). A sequence of positive numbers $\omega(n)\nearrow\infty$ is said to be an a.e. convergence Weyl multiplier (shortly C-multiplier) if every series 
	\begin{equation}\label{a1}
	\sum_{n=1}^\infty a_n\phi_n(x),
	\end{equation}
with coefficients satisfying the condition $\sum_{n=1}^\infty a_n^2\omega(n)<\infty$ is a.e. convergent. If such series converge  unconditionally a.e., then we say $\omega(n)$ is an a.e unconditional convergence Weyl multiplier (UC-multiplier) for $\Phi$.
\end{definition}
Note that Menshov \cite{Men} and Rademacher \cite{Rad} used estimate \e{1} to prove that the sequence $\log^2 n$ is a C-multiplier for any orthonormal system. Likewise, from \trm{T1}, we will deduce the following.
\begin{corollary}\label{C1}
 If $\ZF=\{f_n\}$ is a martingale difference, then $\log n$ is a C-multiplier for any system of $L^2$-normalized \textit{non-overlapping} $\ZF$-polynomials
	\begin{equation}
	p_n(x)=\sum_{j\in G_n}c_jf_j(x),\quad n=1,2,\ldots,
	\end{equation}
	where $G_n\subset \ZN$ are finite and pairwise disjoint.
\end{corollary}
The following result is interesting and it immediately follows from \cor{C1}. 
\begin{corollary}\label{C3}
	The sequence $\log n$ is C-multiplier for any rearrangement of a generalized Haar system.
\end{corollary}
\begin{corollary}\label{C2}
	Let $\{p_n\}$ be a sequence of $L^2$-normalized non-overlapping polynomials with respect to a martingale difference. If $\omega(n)/\log n$ is increasing and
	\begin{equation}\label{omega}
	\sum_{n=1}^\infty \frac{1}{n\omega(n)}<\infty,
	\end{equation}
then $\omega(n)$ is UC-multiplier for $\{p_n\}$.
\end{corollary}

The optimality of $\log n$ in \cor{C3}  as well as  condition \e{omega} in \cor{C2} both follows from some results of Ulyanov for classical Haar system (see \cite{Uly1},\cite{Uly2} or \cite{KaSa} ch. 2 Theorem 17).  In particular, the paper \cite{Uly1} proves that \e{omega} is a necessary and sufficient condition for a sequence $\omega(n)\nearrow\infty$ to be an UC-multiplier for the classical Haar system. 

We prove \trm{T1} using a good-$\lambda$ inequality due to Chang-Wilson-Wolff  \cite{CWW}. See also \cite{GHS}, where the same method has been first applied in the study of maximal functions of  Mikhlin-H\"{o}rmander multipliers. 

\begin{remark}
	Recall that an orthonormal system $\Phi$ is said to be a convergence system if $\omega(n)\equiv 1$ is a C-multiplier for $\Phi$.	It was proved by Koml\'{o}s-R\'{e}v\'{e}sz \cite{KoRe} that if an orthonormal system $\Phi=\{\phi_n\}\subset L^2(0,1)$ satisfies $\|\phi_n\|_4\le M,\,n=1,2,\ldots,$ and we have
	\begin{equation}\label{a20}
	\int_0^1\phi_{n_1}\phi_{n_2}\phi_{n_3}\phi_{n_4}=0
	\end{equation}
	for any choice of different indexes $n_1,n_2,n_3,n_4$, then $\Phi$ is a convergence system. One can check that systems of non-overlapping martingale difference polynomials satisfy \e{a20}. Thus, with the extra condition $\|p_n\|_4\le M$ in \cor{C1} we can claim that $\{p_n\}$ is a convergence system.
\end{remark}
\begin{question}
	Is the additional condition  $\|p_n\|_p\le M$ in \cor{C1}, with a fixed $2<p<4$, is sufficient for $\{p_n\}$ to be a convergence system?
\end{question}

\section{Measure-preserving transformations}
 A mapping $\tau:[0,1)\to[0,1)$ is said to be measure-preserving (MP) transformation if $|\tau^{-1}(A)|=|A|$ for any Lebesgue measurable set $A\subset [0,1)$. A set in $[0,1)$ is said to be simple, if it is a finite union of intervals (of the form $[\alpha,\beta)$).  Let $a$ be a simple set. One can easily check, that the function
\begin{equation}
\xi_a(x)=\frac{|[0,x)\cap a|}{|a|}
\end{equation}
defines a one to one mapping from $a$ to $[0,1)$, such that $|\xi_a(E)|=|E|/|a|$ for any Lebesgue measurable set $E\subset a$. Given integer $n\ge 1$ the mapping $\eta_n(x)=\{nx\}$ defines an MP-transformation of $[0,1)$. Observe that if $a$ is a simple set, then for any integer $n\ge 1$ the mapping  
\begin{equation}
u_{a, n}(x)=\left\{\begin{array}{lcl}
((\xi_a)^{-1}\circ \eta_n \circ \xi_a)(x)&\hbox{ if }& x\in a,\\
x&\hbox{ if }&x\in [0,1)\setminus a ,
\end{array}
\right.
\end{equation}
determines an MP-transformation of $[0,1)$ that maps the set $a$ to itself. Moreover, for any functions $f,g\in L^2(0,1)$ we have
\begin{equation}\label{b5}
\lim_{n\to\infty}\int_a f(u_{a ,n}(x))g(x)dx=\int_a f(x)dx\cdot \int_a g(x)dx
\end{equation}
that is a well-known standard argument. A partition of $[0,1)$ is a sequence of pairwise disjoint sets $\ZA=\{E_k\}$ such that $\cup_kE_k=[0,1)$. We say $\ZA$ is a simple partition if each $E_k$ is simple. Let $\ZA=\{a_j\}$ be a simple partition of $[0,1)$. Given integer $n\ge1$ we consider the MP-transformation
\begin{equation}
u_{\ZA, n}(x)=\sum_ju_{a_j,n}(x)\cdot \ZI_{a_j}(x)
\end{equation}
that maps every $a_j$ to itself. This is an MP-transformation on $[0,1)$ that maps each set $a_j$ to itself and from \e{b5} it follows that
\begin{equation}\label{b7}
\lim_{n\to\infty}\int_0^1f(u_{\ZA,n}(x))g(x)dx=\sum_j\int_{a_j}f(x)dx\cdot \int_{a_j}g(x)dx
\end{equation}
for any functions $f,g\in L^2(0,1)$. An MP-transformation $\tau$ is said to be simple if $\tau^{-1}(a)$ is simple set whenever $a$ is simple. Obviously all above described MP-transformations are simple.

A sequence $\ZA_n$, $n=1,2,\ldots$, of partitions of $[0,1)$ is said to be a filtration if 
any $A\in \ZA_n$  is a union of some sets from $\ZA_{n+1}$ called children of $A$. A martingale difference based on a filtration $\{\ZA_n:\, n=1,2,\ldots \}$ is a sequence of functions $f_n\in L^1(0,1)$, satisfying the conditions 
\begin{enumerate}
	\item Every function $f_n$ is constant on each $A\in \ZA_n$.
	\item We have $\int_Af_n=0$ for any $A\in \ZA_{n-1}$, $n\ge 2$.
\end{enumerate}
Consider a filtration $\{\ZA_n\}$ for which 1) $\ZA_1$ consists of a single element $[0,1)$, 2) each element $A\in \ZA_n$ has only two children intervals in $\ZA_{n+1}$, 3) $\max_{A\in \ZA_n}|A|\to 0$ as $n\to\infty$. A generalized Haar system is a $L^2$-normalized martingale difference based on such filtration. If two children intervals of any $A\in \ZA_n$ are equal, then it gives a signed classical Haar system. It is well-known that any generalized Haar system is complete.

We say that a function system $\{\tilde f_n\}$ is a transformation of another system $\{f_n\}$ if for every choice of numbers $m_k\in \ZN$ and $\lambda_k\in \ZR$ it holds the equality
 \begin{equation}\label{b8}
|\{f_{m_k}(x)>\lambda_k, \, k=1,2,\ldots,n \}|=|\{\tilde f_{m_k}(x)>\lambda_k, \, k=1,2,\ldots,n \}|.
\end{equation}
For example, this relation occurs when $\tilde f_k(x)= f_k(\tau(x))$ for some MP-transformation $\tau$.

The following lemma is an extension of a lemma of Olevksii \cite{Ole} (see also \cite{KaSa}, ch. 10, Lemma 1) proving the same for the classical Haar system. 
\begin{lemma}\label{L1}
	Let $\Phi=\{\phi_k(x)\}$ be a complete orthonormal system and $\ZF=\{f_n\}$ be a martingale difference based on a filtration consistiong of intervals. Then for any sequence of numbers $\varepsilon_k>0$ there exists a transformation $\tilde \ZF=\{\tilde f_n\}$ of the system $\ZF$ and a sequence of non-overlapping $\Phi$-polynomials $p_k$ such that 
	\begin{equation}\label{b3}
	\|\tilde f_k-p_k\|_2<\varepsilon_k,\quad k=1,2,\ldots.
	\end{equation}
\end{lemma}
\begin{proof}
	We shall realize the constructions of sequences $\tilde f_k$ and $p_k$ by induction. First, we take $\tilde f_1=f_1$. Approximation of $f_1$ by a $\Phi$-polynomial $p_1$ gives \e{b3} for $k=1$ that is the base of induction. Then suppose that we have already defined $\tilde f_k, \,p_k$, $k=1,2,\ldots, l$, satisfying the condition \e{b3} such that $\tilde f_k(x)=f_k(\tau_l(x))$, $k=1,2,\ldots,l$,  where $\tau_l$ is a for a simple MP-transformation (maps a simple set to a simple set) . 
	Let $\ZA=\{a_j\}$ be the partition of $[0,1)$ that is formed by the maximal sets, where each function $\tilde f_k$, $k=1,2,\ldots, l$ is constant. Clearly each $a_j$ is a simple set. Since $u_{\ZA, n}$ maps each $a_j$ to itself,  $\tau_{l+1}=\tau_l\circ u_{\ZA, n}$ determines a simple MP-transformation so that $f_k(	\tau_{l+1}(x))=f_k(	\tau_{l}(x))=\tilde f_k(x)$, $k=1,2,\ldots, l$, and
	\begin{equation}\label{b9}
	\int_{\alpha_i}f_{l+1}(\tau_{l}(x))dx=0, \quad i=1,2,\ldots.
	\end{equation}
	From \e{b7} and \e{b9} it follows that
	\begin{align}\label{b6}
	\lim_{n\to\infty}\int_0^1f_{l+1}(	\tau_{l+1}(x))\phi_i(x)dx&= \lim_{n\to\infty}\int_0^1f_{l+1}(\tau_l\circ u_{\alpha, n})(x)\phi_i(x)dx\\
	&=\sum_i\int_{\alpha_i}f_{l+1}(\tau_{l}(x))dx\int_{\alpha_i}\phi_j(x)dxdx=0
	\end{align}
	for any $i=1,2,\ldots $.  We will chose $n$ bigger enough and define $\tilde f_{l+1}(x)=f_{l+1}(\tau_{l+1}(x))$. Let $c_i$ be the Fourier coefficients of the function $\tilde f_{l+1}$ in system $\Phi$. Suppose that each polynomial $p_k$, $k=1,2,\ldots,l$, is a linear combination of functions $\phi_j$ , $j=1,2,\ldots,m$. From \e{b6} it follows that for a bigger enough $n$ we have $\sum_{i=1}^mc_i^2<\varepsilon_{l+1}^2/4$.
	Then we can chose an integer $r>m$ such that $\sum_{i=r+1}^\infty c_i^2<\varepsilon_{l+1}^2/4$.	Define
	\begin{equation}
	p_{l+1}(x)=\sum_{i=m+1}^rc_i\phi_i(x).
	\end{equation}
	Since $\Phi$ is a complete system, one can easily check that \e{b3} is satisfied for $k=l+1$ that finalizes the induction and so the proof of lemma.
\end{proof}

\section{Proof of \trm{T1}}

We will first prove the theorem for the classical Haar system. Let $h_n$ be the $L^2$-normalized classical Haar system. For a given function $f\in L^1(0,1)$ let
$\sum_{k=1}^\infty a_kh_k$ be the Fourier-Haar series of $f$. Recall the maximal and the square functions operators defined by
\begin{equation*}
\ZM f(x)=\sup_{n\ge 1}\left|\sum_{k=1}^na_kh_k(x)\right|,\quad Sf(x)=\left(\sum_{k=1}^\infty a_k^2h_k^2(x)\right)^{1/2}.
\end{equation*}
It is well known the boundedness of both operators on $L^p$, $1<p<\infty$. A key point in the proof of \trm{T1} is the following good-$\lambda$ inequality due to Chang-Wilson-Wolff  (see \cite{CWW}, Corollary 3.1):
\begin{multline}\label{CWW}
|\{x\in[0,1):\,\ZM f(x)>\lambda,\, Sf(x)<\varepsilon\lambda\}|\\\lesssim\exp\left(-\frac{c}{\varepsilon^2}\right)|\{\ZM f(x)>\lambda/2\}|,\,\lambda>0,\,0<\varepsilon <1.
\end{multline}
So let $p_k$, $k=1,2,\ldots,n$, be a monotonic sequence of Haar polynomials. We have
$|g(x)|\le \ZM g(x)$ a.e. for any function $g\in L^1$, as well as $Sp_k(x)\le Sp_n(x)$, $k=1,2,\ldots,n$. Thus, applying inequality \e{CWW} with $\varepsilon_n=(c/\ln n)^{1/2}$, we obtain 
\begin{align}\label{b2}
|\{|p_k(x)|>&\lambda,\, Sp_n(x)\le\varepsilon_n\lambda\}|\\
&\lesssim\exp\left(-\frac{c}{\varepsilon_n^2}\right)|\{\ZM p_k(x)>\lambda/2\}|.
\end{align}
For $p^*(x)=\max_{1\le m\le n}|p_m(x)|$ we obviously have
\begin{align}
\{ p^*(x)>\lambda\}&\subset \{p^*(x)>\lambda,\, Sp_n(x)\le \varepsilon_n\lambda\}\\
&\cup  \{Sp_n(x)> \varepsilon_n\lambda\}=A(\lambda)\cup B(\lambda),
\end{align}
and thus
\begin{equation}
\|p^*\|_2^2\le 2\int_0^\infty\lambda |A(\lambda)|d\lambda+2\int_0^\infty\lambda |B(\lambda)|d\lambda.
\end{equation}
From \e{b2} it follows that
\begin{align}
\int_0^\infty\lambda|A(\lambda)|d\lambda&\le \sum_{m=1}^n\int_0^\infty\lambda |\{ |p_m|>\lambda,\, Sp_n\le \varepsilon_n\lambda\}|d\lambda\\
&\le \exp\left(-\frac{c}{\varepsilon_n^2}\right)\sum_{m=1}^n\int_0^\infty\lambda |\{ \ZM p_m>\lambda/2\}|d\lambda\\
&\lesssim \frac{1}{n}\sum_{m=1}^n\|\ZM p_m\|_2^2\\
&\lesssim \frac{1}{n}\sum_{m=1}^n\| p_m\|_2^2\\
&\le \|p_n\|_2^2.
\end{align}
Combining this and 
\begin{align*}
2\int_0^\infty\lambda|B(\lambda)|d\lambda&=\varepsilon_n^{-2}\|Sp_n\|_2^2\lesssim \log n \cdot \|p_n\|_2^2,
\end{align*}
we get
\begin{equation}\label{b4}
\|p^*\|_2=\left\|\max_{1\le m\le n}|p_m(x)|\right\|_2\lesssim \sqrt{\log n}\cdot \|p_n\|_2
\end{equation} 
that proves the theorem for the Haar system. Clearly we will have the same bound also for any transformation of the Haar system. To proceed the general case we suppose that $\ZF=\{f_n\}$ is an arbitrary martingale difference and let  
\begin{equation}
F_k=\sum_{j\in G_k}c_jf_j,\quad k=1,2,\ldots, n,
\end{equation}
be an arbitrary monotonic sequence of $\ZF$-polynomials. Apply \lem{L1}, choosing $\Phi $ to be the Haar classical system and $\varepsilon_j=\varepsilon$ for $j\in G_n$. So we get \e{b3} for non-overlapping Haar polynomials $p_k$.  Denote $\tilde F_k=\sum_{j\in G_k}c_j\tilde f_j$. Obviously,
\begin{equation}
P_k=\sum_{j\in G_k}c_jp_j,\quad k=1,2,\ldots, n,
\end{equation}
forms a monotonic sequence of Haar polynomials. For a small enough $\varepsilon$ we will have
\begin{align}
\|\tilde F_k-P_k\|_2&\le\left(\sum_{j\in G_n}c_j^2 \right)^{1/2}\left(\sum_{j\in G_n}\varepsilon_j^2 \right)^{1/2}\\
&=\varepsilon \sqrt{\# (G_n)}\left(\sum_{j\in G_n}c_j^2 \right)^{1/2}\le \frac{\|P_n\|_2}{n}.
\end{align}
Therefore, taking into account that the theorem is true for the Haar system, we get
\begin{align}
	\left\|\max_{1\le m\le n}\left|F_k\right|\,\right\|_2&=	\left\|\max_{1\le m\le n}\left|\tilde F_k\right|\,\right\|_2\le \left\|\max_{1\le m\le n}\left|P_k\right|\,\right\|_2+\|P_n\|_2\\
	&\lesssim \sqrt{\log n}\cdot \|P_n\|_2\lesssim\sqrt{\log n}\cdot \|F_n\|_2.
\end{align}
This completes the proof of theorem.

\section{Proof of \trm{T2}}

The upper bound $\kappa_n(\ZH)\lesssim \sqrt{\log n}$ follows from \trm{T1}. The lower bound
	\begin{equation}\label{b1}
	\kappa_n(\ZH)\gtrsim \sqrt{\log n}
	\end{equation}
	for the classical Haar system follows from the Nikishin-Ulyanov \cite{NiUl} inequality 
	\begin{equation*}
	\left\|\sup_{1\le m\le n}\left|\sum_{k=1}^ma_k\chi_{\sigma(k)}\right|\,\right\|_2\gtrsim \sqrt{\log n}\cdot \left(\sum_{k=1}^na_k^2\right)^{1/2},
	\end{equation*}
	valid for appropriate coefficients $a_k$ and permutation $\sigma$ of the numbers $\{1,2,\ldots,n\}$. We will have the same estimate \e{b1} also for any transformation of the classical Haar system.
	Then we apply Olevskii lemma (\cite{KaSa}, ch. 10, Lemma 1), that is the case of \lem{L1} when $\ZF$ coincides with the classical Haar system. So we get a transformed Haar system $\{\tilde h_n\}$ 
	and a sequence of non-overlapping $\Phi$-polynomials $p_k$ such that 
	\begin{equation}
	\|\tilde h_k-p_k\|_2<\varepsilon_k,\quad k=1,2,\ldots.
	\end{equation}
	Since $\varepsilon_k$'s here can be arbitrarily small,  one can conclude $\kappa_n(\Phi)\ge \kappa_n(\ZH)$. Combining this and \e{b1} we get the following.
	\begin{proposition}\label{T4}
		If $\Phi$ is a complete orthonormal system, then $\kappa_n(\Phi)\gtrsim \sqrt{\log n}$.
	\end{proposition}
	 Since any generalized Haar system is complete, the lower bound \e{b1} immediately follows from \pro{T4}.
\section{Proof of corollaries}
\begin{lemma}[\cite{KaSt}, Theorem 5.3.2]\label{L2}
 Let $\{\phi_n(x)\}$ be an orthonormal system and $\omega(n)\nearrow\infty $ be a sequence of positive numbers. If an increasing sequence of indexes $n_k$ satisfy the bound $\omega(n_k)\ge k$, then the condition $\sum_{k=1}^\infty a_k^2\omega(k)<\infty$ implies a.e. convergence of sums $ \sum_{j=1}^{n_k}a_j\phi_j(x)$
 as $k\to\infty$.
\end{lemma}

\begin{proof}[Proof of \cor{C1}]
	Consider the series 
	\begin{equation}\label{d1}
	\sum_{k=1}^\infty a_kp_k(x)
	\end{equation}
	with coefficients satisfying the condition $\sum_{k=1}^\infty a_k^2\log k<\infty$ and denote $S_n=\sum_{k=1}^np_k$.
	Since $\omega(n)=\log n$ satisfies the condition $\omega(2^k)\ge k$, from \lem{L2} we have  a.e. convergence of subsequences $S_{2^k}(x)$. So we just need to show that \begin{equation}\label{d2}
	\delta_k(x)=\max_{2^k< n\le 2^{k+1}}|S_n(x)-S_{2^k}(x)|\to 0\text { a.e. as }k\to\infty.
	\end{equation}
	We have 
\begin{equation*}
	\|\delta_k\|_2\le \kappa_{2^k}(\ZF)\left(\sum_{j=2^k+1}^{2^{k+1}}a_j^2\right)^{1/2}\lesssim \sqrt{k}\left(\sum_{j=2^k+1}^{2^{k+1}}a_j^2\right)^{1/2}. 
\end{equation*}
So we get
	\begin{equation*}
	\sum_{k=1}^\infty\|\delta_k\|_2^2\le \sum_{k=1}^\infty k\sum_{j=2^k+1}^{2^{k+1}}a_j^2\le \sum_{j=1}^\infty a_j^2\log j<\infty,
	\end{equation*}
	which implies \e{d2}.
\end{proof}
To prove the next corollary we will need another lemma.
\begin{lemma}[\cite{Uly3}, \cite{Pol}]\label{L3}
	Let $u(n)$ be a C-multiplier for any rearrangement of the orthonormal system $\Phi=\{\phi_n(x)\}$. If an increasing sequence of positive numbers $\delta(k)$ satisfies the condition
	\begin{equation}\label{d3}
	\sum_{k=1}^\infty\frac{1}{\delta(k)k\log k}<\infty,
	\end{equation}
	then $\delta(n)u(n)$ turnes to be a UC-multiplier for $\Phi$.
\end{lemma}
\begin{proof}[Proof of \cor{C2}]
	According to \cor{C1} $u(n)=\log n$ is a C-multiplier for the systems of non-overlapping MD-polynomials and their rearrangements. By the hypothesis of \cor{C2} the sequence $\delta(n)=\omega(n)/\log n$ is increasing and satisfies \e{d3}. Thus, the combination of \cor{C1} and \lem{L3} completes the proof.
\end{proof}

\bibliographystyle{plain}


\end{document}